\RequirePackage{fix-cm} 
\documentclass[a4paper, twoside, reqno, dvips, 12pt]{amsart}
\usepackage{fixltx2e}   

\usepackage{etex}

\usepackage[latin1]{inputenc}
\usepackage[T1]{fontenc}

\usepackage{eucal}
\usepackage{esint}
\usepackage{dsfont}
\usepackage{xspace}
\usepackage{amsgen}
\usepackage{amsthm}
\usepackage{amssymb}
\usepackage{amsmath}
\usepackage{upgreek}
\usepackage{wasysym}
\usepackage{amsfonts}
\usepackage{stmaryrd}
\usepackage{mathtools}

\usepackage{mathrsfs}
\DeclareMathAlphabet{\mathscrbf}{OMS}{mdugm}{b}{n}

\usepackage{a4wide}

\usepackage[dvipsnames, table]{xcolor}
\definecolor{bckg}{RGB}{20.8, 20.8, 20.8}
\definecolor{oneblue}{rgb}{0.0, 0.0, 0.85}
\definecolor{Lightblue}{RGB}{214, 214, 214}
\definecolor{bluepigment}{rgb}{0.2, 0.2, 0.6}
\definecolor{charcoal}{rgb}{0.21, 0.27, 0.31}
\definecolor{denimblue}{rgb}{0.08, 0.38, 0.74}
\definecolor{Lightgray}{rgb}{0.89, 0.89, 0.89}
\definecolor{darkgrey}{rgb}{0.273, 0.281, 0.30}
\definecolor{darkelectricblue}{rgb}{0.33, 0.41, 0.47}

\usepackage{psfrag}
\usepackage{graphicx}
\usepackage{subfigure}
\usepackage{morefloats}
\usepackage{indentfirst}

\usepackage{acronym}
\usepackage{microtype}
\usepackage[labelsep=period,%
            labelfont={bf,sf,color=bluepigment},%
            justification=raggedright]{caption}

\usepackage[perpage, symbol]{footmisc}

\usepackage[usenames, dvipsnames]{pstricks}
\usepackage{epsfig}
\usepackage{pst-grad} 
\usepackage{pst-plot} 

\usepackage[colorlinks,
           urlcolor=oneblue,
           linkcolor=denimblue,
           citecolor=NavyBlue,
           bookmarksopen=false,
           pdfpagemode=UseNone,
           pagebackref]{hyperref}

\usepackage[sort&compress, comma, square, numbers]{natbib}

\usepackage[explicit]{titlesec}

\titleformat{\section}
  {\color{NavyBlue}\Large\sffamily\bfseries}
  {}
  {0em}
  {\colorbox{bckg!5}{\parbox{\dimexpr\linewidth-2\fboxsep\relax}{\centering\thesection. #1}}}
  [\vspace*{0.33em}]

\titleformat{name=\section,numberless}
  {\color{NavyBlue}\Large\sffamily\bfseries}
  {}
  {0.0em}
  {\colorbox{bckg!10}{\parbox{\dimexpr\linewidth-2\fboxsep\relax}{\centering#1}}}
  [\vspace*{0.33em}]

\titleformat{\subsection}
  {\color{NavyBlue}\large\sffamily\bfseries}
  {}
  {0.0em}
  {\colorbox{bckg!5}{\parbox{\dimexpr\linewidth-2\fboxsep\relax}{\centering\thesubsection. #1}}}
  [\vspace*{0.33em}]

\titleformat{name=\subsection,numberless}
  {\color{NavyBlue}\Large\sffamily\bfseries}
  {}
  {0em}
  {\colorbox{bckg!10}{\parbox{\dimexpr\linewidth-2\fboxsep\relax}{\centering#1}}}
  [\vspace*{0.33em}]

\titleformat{\subsubsection}
  {\color{bluepigment}\sffamily\normalsize\bfseries}
  {\thesubsubsection}
  {0.5em}
  {#1}
  [\vspace*{0.33em}]

\titleformat{\paragraph}[runin]
  {\color{bluepigment}\sffamily\small\bfseries}
  {}
  {0em}
  {#1}

\titlespacing{\section}{1.0em}{1.5em plus 2pt minus 2pt}%
{1.0em plus 2pt minus 2pt}[0em]
\titlespacing{\subsection}{1.0em}{1.5em plus 2pt minus 2pt}%
{1.0em}[0em]
\titlespacing{\subsubsection}{1.0em}{1.5em plus 2pt minus 2pt}%
{1.0em plus 2pt minus 2pt}[0em]

\usepackage{titletoc}

\contentsmargin{0.5em}
\newlength{\tocsep} 
\titlecontents{section}[\tocsep]
  {\addvspace{10pt}\bfseries\sffamily}
  {\contentslabel[\thecontentslabel]{\tocsep}}
  {}
  {\ \titlerule*[0.75pc]{.}\ \thecontentspage}
  []
\titlecontents{subsection}[\tocsep]
  {\addvspace{8pt}\sffamily}
  {\contentslabel[\thecontentslabel]{\tocsep}}
  {}
  {\ \titlerule*[0.5pc]{.}\ \thecontentspage}
  []
\titlecontents*{subsubsection}[\tocsep]
  {\addvspace{2pt}\footnotesize\sffamily}
  {}
  {}
  {\ \titlerule*[0.35pc]{.}\ \thecontentspage}
  [\\*]

\makeatletter
\def\@setauthors{%
  \begingroup
  \def\thanks{\protect\thanks@warning}%
  \trivlist
  \centering\footnotesize \@topsep30\p@\relax
  \advance\@topsep by -\baselineskip
  \item\relax
  \author@andify\authors
  \def\\{\protect\linebreak}%
  \textsc{\normalsize\textcolor{darkelectricblue}{\authors}}%
  \ifx\@empty\contribs
  \else
    ,\penalty-3 \space \@setcontribs
    \@closetoccontribs
  \fi
  \endtrivlist
  \endgroup
}
\def\@settitle{\begin{center}%
  \baselineskip14\p@\relax
    \bfseries
    \textsc{\Large\textcolor{charcoal}{\@title}}
  \end{center}%
}
\makeatother

\usepackage{enumitem}
\setlist[description]{%
  topsep=30pt,               
  itemsep=5pt,               
  font={\bfseries\sffamily\color{NavyBlue}}, 
}

\usepackage{fancyhdr}
\usepackage{lastpage}

\newcommand*\Title{\textcolor{bluepigment}{How to overcome the CFL?}}
\newcommand*\Authors{\textcolor{bluepigment}{D.~Dutykh}}
\newcommand*{\plogo}{\textcolor{gray}{{\texttt{arXiv.org} / \textsc{hal}}}} 

\numberwithin{equation}{section}

\newtheorem{remark}{Remark}

\renewcommand{\nu}{\upnu}
\renewcommand{\tau}{\uptau}

\newcommand{\R}{\mathds{R}}

\newcommand{\ue}{\mathrm{e}}
\newcommand{\ui}{\mathrm{i}}
\newcommand{\M}{\mathcal{M}}
\renewcommand{\beta}{\upbeta}
\newcommand{\Ll}{\mathscr{L}}
\renewcommand{\L}{\mathcal{L}}
\renewcommand{\O}{\mathcal{O}}
\renewcommand{\leq}{\leqslant}
\renewcommand{\geq}{\geqslant}

\renewcommand{\alpha}{\upalpha}
\renewcommand{\kappa}{\upkappa}
\renewcommand{\gamma}{\upgamma}
\newcommand{\x}{\boldsymbol{x}}
\renewcommand{\omega}{\upomega}
\renewcommand{\lambda}{\uplambda}
\newcommand{\const}{\mathrm{const}}

\renewcommand{\Re}{\operatorname{Re}}
\renewcommand{\Im}{\operatorname{Im}}

\newcommand{\ie}{\emph{i.e.}\xspace}
\newcommand{\eg}{\emph{e.g.}\xspace}
\newcommand{\etc}{\emph{etc.}\xspace}

\newcommand{\scal}{\boldsymbol{\cdot}}
\newcommand{\grad}{\boldsymbol{\nabla}}

\newcommand{\abs}[1]{\lvert\, #1\, \rvert}

\newcommand{\norm}[1]{\lVert\, #1\, \rVert}

\newcommand{\pd}[2]{\frac{\partial #1}{\partial\/ #2}}
\newcommand{\pdd}[2]{\dfrac{\partial #1}{\partial\/ #2}}

\newcommand{\eqdef}{\mathop{\stackrel{\,\mathrm{def}}{:=}\,}}

\newcommand{\half}{{\textstyle{1\over2}}}

\newcommand{\sixth}{{\textstyle{1\over6}}}

\newcommand{\twelwth}{{\textstyle{1\over12}}}
\newcommand{\twothird}{{\textstyle{2\over3}}}

\newcommand{\threefourth}{{\textstyle{3\over4}}}

\usepackage{acronym}

\begin{document}

\title[\Title]{How to overcome the Courant--Friedrichs--Lewy condition of explicit discretizations?}

\author[D.~Dutykh]{Denys Dutykh$^*$}
\address{LAMA, UMR 5127 CNRS, Universit\'e Savoie Mont Blanc, Campus Scientifique, 73376 Le Bourget-du-Lac Cedex, France}
\email{Denys.Dutykh@univ-savoie.fr}
\urladdr{http://www.denys-dutykh.com/}
\thanks{$^*$ Corresponding author}



\begin{titlepage}
\thispagestyle{empty} 
\noindent
{\Large Denys \textsc{Dutykh}}\\
{\it\textcolor{gray}{CNRS--Universit\'e Savoie Mont Blanc, France}}
\\[0.08\textheight]

\vspace*{3cm}

\colorbox{Lightblue}{
  \parbox[t]{1.0\textwidth}{
    \centering\huge\sc
    \vspace*{0.7cm}
    
    \textcolor{bluepigment}{How to overcome the Courant--Friedrichs--Lewy condition of explicit discretizations?}
    
    \vspace*{0.7cm}
  }
}

\vfill 

\raggedleft     
{\large \plogo} 
\end{titlepage}


\newpage
\thispagestyle{empty} 
\par\vspace*{\fill}   
\begin{flushright} 
{\textcolor{denimblue}{\textsc{Last modified:}} \today}
\end{flushright}


\newpage
\maketitle
\thispagestyle{empty}


\begin{abstract}

This manuscript contains some thoughts on the discretization of the classical heat equation. Namely, we discuss the advantages and disadvantages of explicit and implicit schemes. Then, we show how to overcome some disadvantages while preserving some advantages. However, since there is no free lunch, there is a price to pay for any improvement in the numerical scheme. This price will be thoroughly discussed below.

In particular, we like explicit discretizations for the ease of their implementation even for nonlinear problems. Unfortunately, when these schemes are applied to parabolic equations, severe stability limits appear for the time step magnitude making the explicit simulations prohibitively expensive. Implicit schemes remove the stability limit, but each time step requires now the solution of linear (at best) or even nonlinear systems of equations. However, there exists a number of tricks to overcome (or at least to relax) severe stability limitations of explicit schemes without going into the trouble of fully implicit ones. The purpose of this manuscript is just to inform the readers about these alternative techniques to extend the stability limits. It was not written for classical scientific publication purposes.

\bigskip
\noindent \textbf{\keywordsname:} Heat equation; Finite differences; Explicit schemes; Implicit schemes; CFL condition \\


\end{abstract}


\newpage
\tableofcontents
\thispagestyle{empty}


\newpage
\section{Introduction}

In this text we consider the classical linear heat equation:
\begin{equation}\label{eq:heat}
  u_{\,t}\ =\ \nu\,\grad^{\,2}\,u\,,
\end{equation}
where $u(\x,\,t)$ is a quantity being diffused in some domain $\Omega\ \subseteq\ \R^{d}\,$. In physical applications $u(\x,\,t)\,$, $\x\ \in\ \Omega\,$, $t\ >\ 0$ may represent the temperature field, moisture content, vapor concentration, \etc and $\nu\ >\ 0$ is the diffusion coefficient. The subscripts denote partial derivatives, \ie $u_{\,t}\ \eqdef\ \pdd{u(\x,\,t)}{t}\,$. Finally, $\grad^{\,2}\ \equiv\ \grad\scal\grad$ is the classical $d-$dimensional \textsc{Laplace} operator:
\begin{equation*}
  \grad^{\,2}\ \eqdef\ \sum_{i\,=\,1}^{d}\, \pd{{}^{\,2}}{x_{\,i}^{\,2}}\,.
\end{equation*}
The derivation of this equation for \textsc{Brownian} motion process was given by A.~\textsc{Einstein} \cite{Einstein1905} in 1905.

From now on we shall restrict our ambitions on the $1-$dimensional case where $\Omega\ \equiv\ [\,0,\,\ell\,]\ \subseteq\ \R^{1}$ and equation \eqref{eq:heat} correspondingly becomes:
\begin{equation}\label{eq:heat1d}
  u_{\,t}\ =\ \nu\,u_{\,x\,x}\,.
\end{equation}
This equation has to be supplemented by one initial
\begin{equation*}
  u\bigl\vert_{t\,=\,0}\ =\ u_{\,0}\,(x)\,,
\end{equation*}
and two boundary conditions:
\begin{align}\label{eq:bcl}
  \Phi_{\,\mathrm{l}}\,\bigl(t,\,u(t,\,0),\,u_x(t,\,0)\bigr)\ &=\ 0\,, \\
  \Phi_{\,\mathrm{r}}\,\bigl(t,\,u(t,\,\ell),\,u_x(t,\,\ell)\bigr)\ &=\ 0\,.\label{eq:bcr}
\end{align}
The functions $\Phi_{\,\mathrm{l},\,\mathrm{r}}\,(\,\bullet\,)$ have to be specified depending on the practical situation in hands. For example, if we have the \textsc{Dirichlet}-type condition on the left boundary then
\begin{equation*}
  \Phi_{\,\mathrm{l}}\,\bigl(t,\,u(t,\,0),\,u_{\,x}\,(t,\,0)\bigr)\ \equiv\ u(t,\,0)\ -\ u^\circ(t)\ =\ 0\,,
\end{equation*}
where $u^\circ(t)$ is a prescribed function of time. Often, it is assumed that $u^\circ(t)\ \equiv\ \const\,$. The homogeneous \textsc{Neumann}-type condition on the right looks like
\begin{equation*}
  \Phi_{\,\mathrm{r}}\,\bigl(t,\,u(t,\,\ell),\,u_{\,x}\,(t,\,\ell)\bigr)\ \equiv\ u_{\,x}\,(t,\,\ell)\ =\ 0\,.
\end{equation*}


\subsection{Some healthy criticism}

Obviously, the heat equation \eqref{eq:heat} is a simplified model obtained after a series of idealizations and simplifications. As a result, equation \eqref{eq:heat} is linear and its \textsc{Green} function\footnote{The \textsc{Green} function is also known as the fundamental solution. More precisely, it solves the following problem:
\begin{align*}
  G_{\,t}\ &=\ \nu\,G_{\,x\,x}\,, \qquad x\ \in\ \R\,, \\
  G\bigr\vert_{\,t\,=\,0}\ &=\ \delta\,(x)\,,
\end{align*}
where $\delta(x)$ is \textsc{Dirac} delta function.} can be computed analytically\footnote{This expression is false. To be checked later!}
\begin{equation*}
  G(x,\,t)\ =\ \frac{1}{\sqrt{4\,\pi\,\nu\,t}}\;\ue^{-\,\dfrac{x^{\,2}}{4\,\nu\,t}}
\end{equation*}
In particular, one can see that for any sufficiently small $t\ >\ 0$ the function $G(x,\,t)$ is not of \emph{compact support}. In other words, the information about a point source initially localized at $x\ =\ 0$ spreads instantly over the whole domain. Of course, the infinite speed of information propagation is physically forbidden. This non-physical feature of heat equation's solutions is a consequence of simplifying assumptions made during the derivation. We just mention here that some nonlinear versions of the heat equation do have fundamental solutions with compact support.

However, in practice the solutions to Partial Differential Equations (PDEs) such as \eqref{eq:heat} are computed numerically than constructed analytically. That is why we can hope to correct some non-physical features of the heat equation solution at the discrete level. This is the main topic of the present manuscript.


\bigskip
\paragraph*{Manuscript organization.}

Below we shall review some classical numerical schemes for the heat equation \eqref{eq:heat} in one spatial dimension in Section~\ref{sec:class}: 
\begin{itemize}
  \item The explicit scheme in Section~\ref{sec:exp}
  \item The implicit scheme in Section~\ref{sec:imp}
  \item The leap-frog scheme in Section~\ref{sec:leap}
  \item The \textsc{Crank}--\textsc{Nicholson} scheme in Section~\ref{sec:crank}
\end{itemize}
In that Section we justify our preference for explicit schemes in time. However, explicit schemes are known to have an important CFL-type stability restrictions on the time step \cite{Courant1928}. That is why in Section~\ref{sec:tricks} we review some not so widely known schemes which allow to overcome the stability limit while still being explicit in time. In particular, we consider the following alternatives:
\begin{itemize}
  \item \textsc{Dufort}--\textsc{Frankel} method in Section~\ref{sec:dufort}
  \item \textsc{Saulyev} method in Section~\ref{sec:saulyev}
  \item Hyperbolization method in Section~\ref{sec:hyper}
\end{itemize}
Finally, the main conclusions and perspectives of the present study are discussed in Section~\ref{sec:disc}.


\section{Classical numerical schemes}
\label{sec:class}

\begin{figure}
  \centering
  \includegraphics[width=0.99\textwidth]{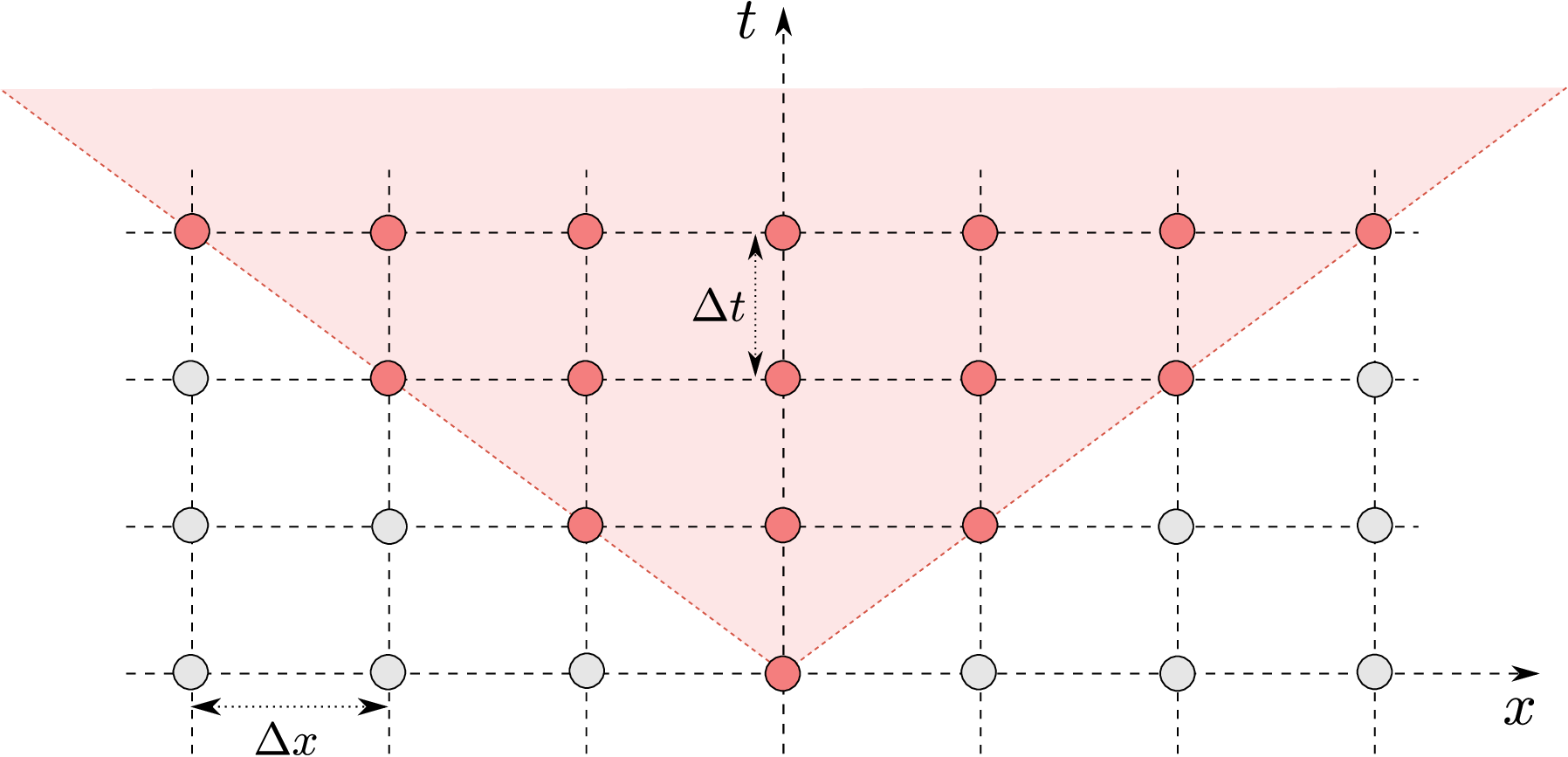}
  \caption{\small\em A schematic representation of the uniform discretization in space and time. Red nodes correspond to non-zero values of the discrete solution. Grey nodes correspond to $u_{\,j}^{\,n}\ \equiv\ 0\,$. The shaded area is an equivalent of the `\emph{light cone}' for the initially activated node.}
  \label{fig:grid}
\end{figure}

In order to describe numerical schemes in simple terms, consider a uniform discretization of the interval $\Omega\ \rightsquigarrow\ \Omega_{\,h}\,$:
\begin{equation*}
  \Omega_{\,h}\ =\ \bigcup_{j\,=\,0}^{N-1} [\,x_{\,j},\;x_{\,j+1}\,]\,, \qquad
  x_{j+1}\ -\ x_{\,j}\ \equiv\ \Delta x\,, \quad \forall j\ \in\ \bigl\{0,\,1,\,\ldots,\,N-1\bigr\}\,.
\end{equation*}
The time layers are uniformly spaced as well $t^{\,n}\ =\ n\,\Delta t\,$, $\Delta t\ =\ \const\ >\ 0\,$, $n\ =\ 0,\,1,\,2,\,\ldots$ The values of function $u(x,\,t)$ in discrete nodes will be denoted by $u_{\,j}^{\,n}\ \eqdef\ u\,(x_{\,j},\,t^{\,n}\,)\,$. The space-time grid is schematically depicted in Figure~\ref{fig:grid}.

\subsection{The Explicit scheme}
\label{sec:exp}

\begin{figure}
  \centering
  \includegraphics[width=0.59\textwidth]{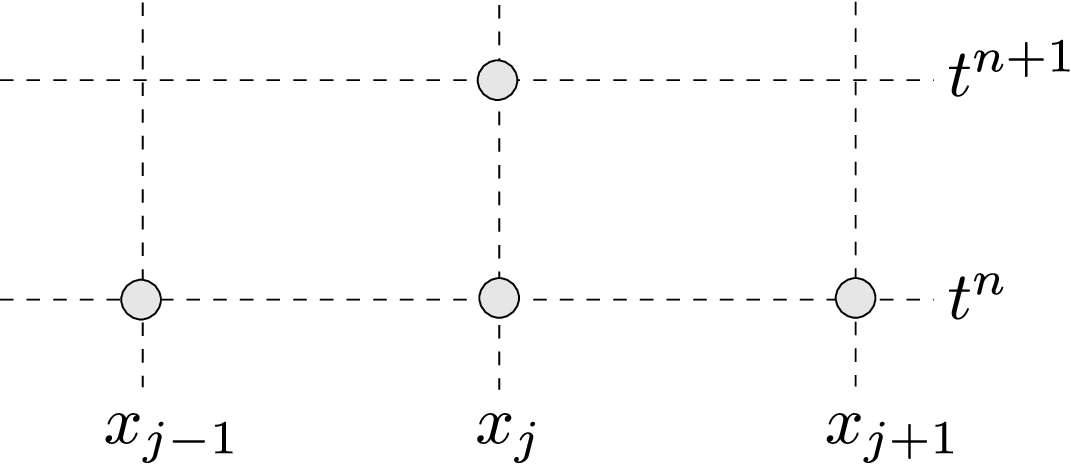}
  \caption{\small\em Stencil of the explicit finite difference scheme \eqref{eq:exp}.}
  \label{fig:exp}
\end{figure}

The standard explicit scheme for the linear heat equation \eqref{eq:heat1d} can be written as
\begin{equation}\label{eq:exp}
  \frac{u_{\,j}^{\,n+1}\ -\ u_{\,j}^{\,n}}{\Delta t}\ =\ \nu\;\frac{u_{\,j-1}^{\,n}\ -\ 2\,u_{\,j}^{\,n}\ +\ u_{\,j+1}^{\,n}}{\Delta x^{\,2}}\,, \qquad j\ =\ 1,\,\ldots,\,N-1\,, \qquad n\ \geq\ 0\,.
\end{equation}
The stencil of this scheme is depicted in Figure~\ref{fig:exp}. In order to complete this discretization we have to find from boundary conditions \eqref{eq:bcl}, \eqref{eq:bcr} the boundary values:
\begin{equation}\label{eq:bc}
  u_{\,0}^{\,n+1}\ =\ \psi_{\,\mathrm{l}}\,(t^{\,n+1},\,u_{\,1}^{\,n+1},\,\ldots\,)\,, \qquad
  u_{\,N}^{\,n+1}\ =\ \psi_{\,\mathrm{r}}\,(t^{\,n+1},\,u_{\,N-1}^{\,n+1},\,\ldots\,)\,,
\end{equation}
where functions $\psi_{\,\mathrm{l},\,\mathrm{r}}\,(\,\bullet\,)$ may depend on adjacent values of the solution whose number depends on the approximation order of the scheme (here we use the second order in space). For example, if the temperature is prescribed on the right boundary, then we simply have
\begin{equation*}
  u_{\,N}^{\,n+1}\ =\ \psi_{\,\mathrm{r}}\,(t^{\,n+1})\ \equiv\ \phi_{\,\mathrm{r}}\,(t^{\,n+1})\,,
\end{equation*}
where $\phi_{\,\mathrm{r}}\,(t)$ is a given function of time. On the other hand, if the heat flux is prescribed on the left boundary $\nu\;\pd{u}{x}\ =\ \phi_{\,\mathrm{l}}\,(t)$ then it can be discretized as
\begin{equation*}
  \nu\;\frac{-3\,u_{\,0}^{\,n+1}\ +\ 4\,u_{\,1}^{\,n+1}\ -\ u_{\,2}^{\,n+1}}{2\,\Delta x}\ =\ \phi_{\,\mathrm{l}}\,(t^{\,n+1})\,.
\end{equation*}

By solving \eqref{eq:exp} with respect to $u_{\,j}^{\,n+1}$ we obtain a discrete dynamical system
\begin{equation*}
  u_{\,j}^{\,n+1}\ =\ u_{\,j}^{\,n}\ +\ \nu\;\frac{\Delta t}{\Delta x^{\,2}}\;\bigl(u_{\,j-1}^{\,n}\ -\ 2\,u_{\,j}^{\,n}\ +\ u_{\,j+1}^{\,n}\bigr)\,,
\end{equation*}
whose starting value is directly obtained from the initial condition:
\begin{equation*}
  u_{\,j}^{\,0}\ =\ u_{\,0}(x_{\,j})\,, \qquad j\ =\ 0,\,1,\,\ldots,\,N\,.
\end{equation*}
It is well-known that scheme \eqref{eq:exp} approximates the continuous operator to order $\O(\Delta t\ +\ \Delta x^{\,2})\,$. The explicit scheme is conditionally stable under the following CFL-type condition:
\begin{equation}\label{eq:cfl}
  \Delta t\ \leq\ \frac{1}{2\,\nu}\;\Delta x^{\,2}\,.
\end{equation}
Unfortunately, this condition is too restrictive for sufficiently fine discretizations.

\subsection{The Implicit scheme}
\label{sec:imp}

\begin{figure}
  \centering
  \includegraphics[width=0.59\textwidth]{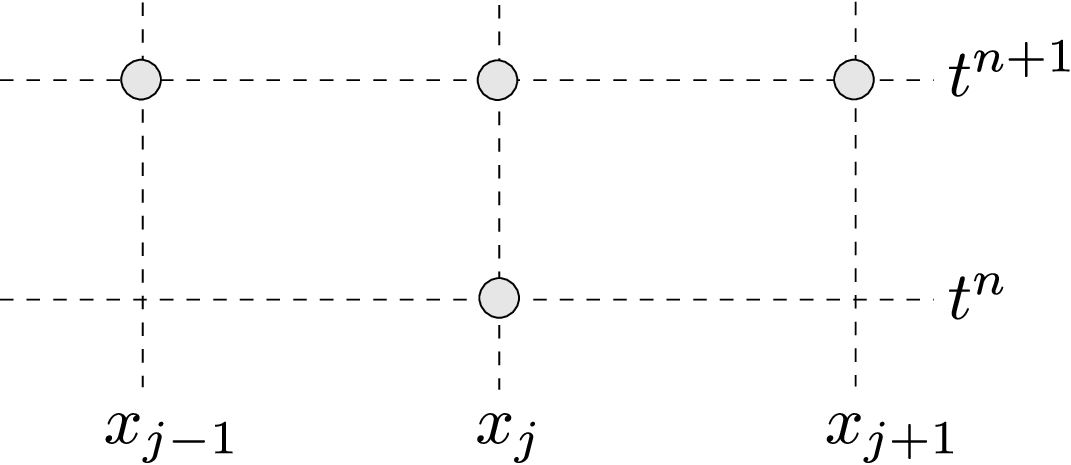}
  \caption{\small\em Stencil of the implicit finite difference scheme \eqref{eq:imp}.}
  \label{fig:imp}
\end{figure}

The implicit scheme for the 1D heat equation \eqref{eq:heat1d} is given by the following relations:
\begin{equation}\label{eq:imp}
  \frac{u_{\,j}^{\,n+1}\ -\ u_{\,j}^{\,n}}{\Delta t}\ =\ \nu\;\frac{u_{\,j-1}^{\,n+1}\ -\ 2\,u_{\,j}^{\,n+1}\ +\ u_{\,j+1}^{\,n+1}}{\Delta x^{\,2}}\,, \qquad j\ =\ 1,\,\ldots,\,N-1\,, \qquad n\ \geq\ 0\,.
\end{equation}
The finite difference stencil of this scheme is depicted in Figure~\ref{fig:imp}. These relations have to be properly initialized and supplemented with numerical boundary conditions \eqref{eq:bc}. In the following Sections we shall not return to the question of initial and boundary conditions in order to focus on the discretization. The scheme \eqref{eq:imp} has the same order of accuracy as the explicit scheme \eqref{eq:exp}, \ie $\O(\Delta t\ +\ \Delta x^{\,2})\,$. However, the implicit scheme \eqref{eq:imp} is unconditionally stable, which constitutes its major advantage. It could be interesting to have also the second order in time as well. This issue will be addressed in the following Sections.

The most important difference with the explicit scheme \eqref{eq:exp} is that we have to solve a tridiagonal system of linear algebraic equations to determine the numerical solution values $\bigl\{u_{\,j}^{\,n+1}\bigr\}_{j\,=\,0}^{N}$ on the following time layer $t\ =\ t^{\,n+1}\,$. It determines the algorithm complexity --- a tridiagonal system of equations can be solved in $\O(N)$ operations (using the simple \textsc{Thomas} algorithm, for example) and it has to be done at every time step.

\subsection{The Leap-frog scheme}
\label{sec:leap}

The leap-frog scheme\footnote{This scheme is called in French as `\emph{le sch\'ema saute-mouton}'.} is obtained by replacing in \eqref{eq:exp} the forward difference in time by the symmetric one, \ie
\begin{equation}\label{eq:leap}
  \frac{u_{\,j}^{\,n+1}\ -\ u_{\,j}^{\,n-1}}{2\,\Delta t}\ =\ \nu\;\frac{u_{\,j-1}^{\,n}\ -\ 2\,u_{\,j}^{\,n}\ +\ u_{\,j+1}^{\,n}}{\Delta x^{\,2}}\,, \qquad j\ =\ 1,\,\ldots,\,N-1\,, \qquad n\ \geq\ 0\,.
\end{equation}
This scheme is second order accurate in space and in time, \ie $\O(\Delta t^{\,2}\ +\ \Delta x^{\,2})\,$. Unfortunately, the leap-frog scheme is unconditionally unstable. It makes it un-exploitable in practice. However, we shall use some modifications of this scheme below.

\begin{figure}
  \centering
  \includegraphics[width=0.59\textwidth]{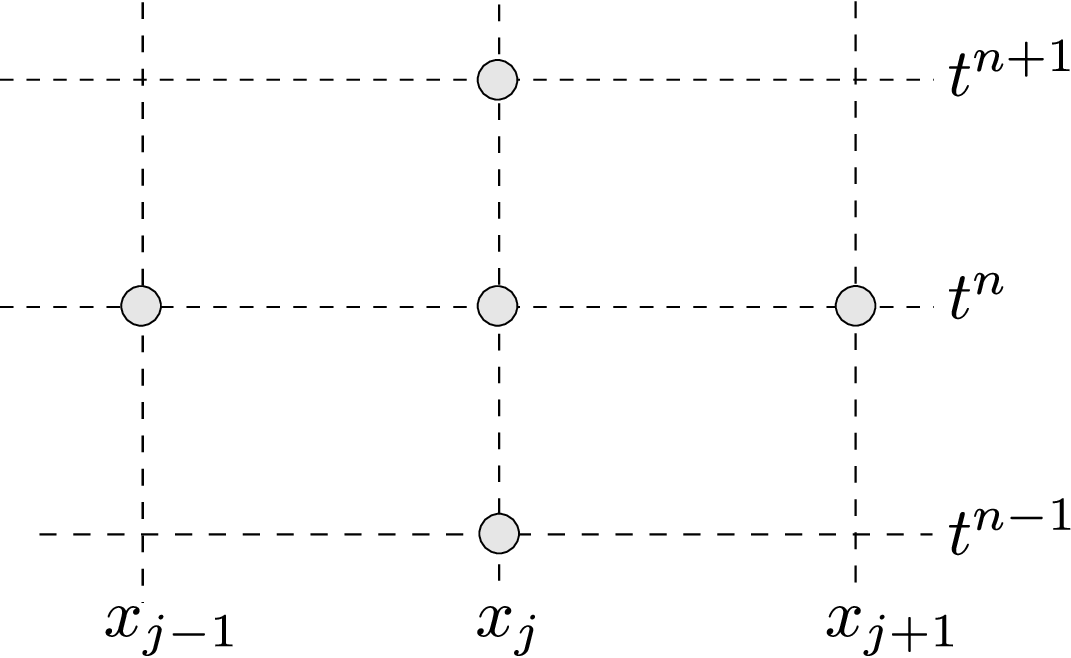}
  \caption{\small\em Stencil of the leap-frog \eqref{eq:leap} and hyperbolic \eqref{eq:hyps} finite difference schemes.}
  \label{fig:leap}
\end{figure}


\subsection{The Crank--Nicholson scheme}
\label{sec:crank}

\begin{figure}
  \centering
  \includegraphics[width=0.59\textwidth]{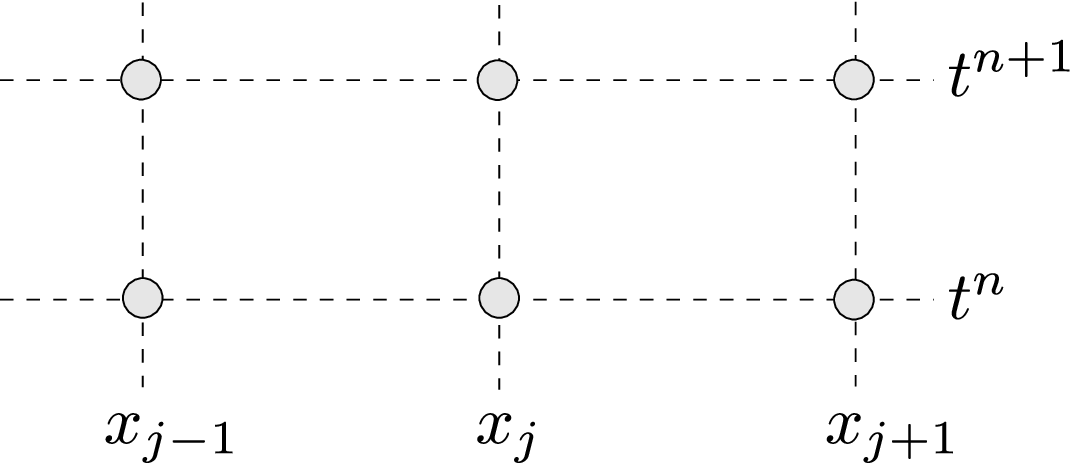}
  \caption{\small\em Stencil of the \textsc{Crank}--\textsc{Nicholson} (CN) finite difference scheme \eqref{eq:cn}.}
  \label{fig:cn}
\end{figure}

We saw above that the first tentative to obtain a scheme with second order accuracy in space \emph{and} in time was unsuccessful (see Section~\ref{sec:leap}). However, a very useful method was proposed by \textsc{Crank} \& \textsc{Nicholson} (CN) and it can be successfully applied to the heat equation \eqref{eq:heat1d} as well:
\begin{multline}\label{eq:cn}
  \frac{u_{\,j}^{\,n+1}\ -\ u_{\,j}^{\,n}}{\Delta t}\ =\ \nu\;\frac{u_{\,j-1}^{\,n}\ -\ 2\,u_{\,j}^{\,n}\ +\ u_{\,j+1}^{\,n}}{2\,\Delta x^{\,2}}\ +\ \nu\;\frac{u_{\,j-1}^{\,n+1}\ -\ 2\,u_{\,j}^{\,n+1}\ +\ u_{\,j+1}^{\,n+1}}{2\,\Delta x^{\,2}}\,, \\ j\ =\ 1,\,\ldots,\,N-1\,, \qquad n\ \geq\ 0\,.
\end{multline}
This scheme is $\O(\Delta t^2\ +\ \Delta x^{\,2})$ accurate and unconditionally stable (similarly to \eqref{eq:imp}). That is why numerical results obtained with the CN scheme will be more accurate than implicit scheme \eqref{eq:imp} predictions. The stencil of this scheme is depicted in Figure~\ref{fig:cn}. The CN scheme has all advantages and disadvantages (except for the order of accuracy in time) of the implicit scheme \eqref{eq:imp}. At every time step one has to use a tridiagonal solver to invert the linear system of equations to determine solution value at the following time layer $t\ =\ t^{\,n+1}\,$.

\subsubsection{Some nonlinear extensions}

Since most of real-world heat conduction models used in building physics are nonlinear, it is worth to discuss some nonlinear extensions of the \textsc{Crank}--\textsc{Nicholson} (CN) scheme. For linear problems CN scheme turns out to be the same as the mid-point and trapezoidal rules for Ordinary Differential Equations (ODEs). Indeed, consider a nonlinear ODE:
\begin{equation}\label{eq:ode}
  \dot{u}\ =\ f(u)\,, \qquad u(0)\ =\ u_{\,0}\,.
\end{equation}
The mid-point and trapezoidal rules consist correspondingly in discretizing \eqref{eq:ode} as follows:
\begin{align*}
  \frac{u^{\,n+1}\ -\ u^{\,n}}{\Delta t}\ &=\ f\,\Bigl(\,\frac{u^{\,n}\ +\ u^{\,n+1}}{2}\,\Bigr)\,, \\
  \frac{u^{\,n+1}\ -\ u^{\,n}}{\Delta t}\ &=\ \frac{f\bigl(u^{\,n}\bigr)\ +\ f\bigl(u^{\,n+1}\bigr)}{2}\,.
\end{align*}
Now, if we set in formulas above $f(u)\ =\ \nu\,\Ll\cdot u\,$, where $\Ll\ \simeq\ \partial_{\,x\,x}$ is the linear operator which represents the second central finite difference, we recover the CN scheme \eqref{eq:cn}.

Consider a non-conservative nonlinear heat equation:
\begin{equation}\label{eq:nheat}
  u_{\,t}\ =\ k(u)\,u_{\,x\,x}\,.
\end{equation}
The straightforward application of the CN scheme to equation \eqref{eq:heat} yields the following scheme:
\begin{multline}\label{eq:cnn}
  \frac{u_{\,j}^{\,n+1}\ -\ u_{\,j}^{\,n}}{\Delta t}\ =\ k(u_{\,j}^{\,n})\;\frac{u_{\,j-1}^{\,n}\ -\ 2\,u_{\,j}^{\,n}\ +\ u_{\,j+1}^{\,n}}{2\,\Delta x^{\,2}}\ +\ k(u_{\,j}^{\,n+1})\;\frac{u_{\,j-1}^{\,n+1}\ -\ 2\,u_{\,j}^{\,n+1}\ +\ u_{\,j+1}^{\,n+1}}{2\,\Delta x^{\,2}}\,, \\ j\ =\ 1,\,\ldots,\,N-1\,, \qquad n\ \geq\ 0\,.
\end{multline}
However, it is less known that one can apply also the cross-\textsc{Crank}--\textsc{Nicholson} (cCN) scheme:
\begin{multline}\label{eq:ccn}
  \frac{u_{\,j}^{\,n+1}\ -\ u_{\,j}^{\,n}}{\Delta t}\ =\ k(u_{\,j}^{\,n+1})\;\frac{u_{\,j-1}^{\,n}\ -\ 2\,u_{\,j}^{\,n}\ +\ u_{\,j+1}^{\,n}}{2\,\Delta x^{\,2}}\ +\ k(u_{\,j}^{\,n})\;\frac{u_{\,j-1}^{\,n+1}\ -\ 2\,u_{\,j}^{\,n+1}\ +\ u_{\,j+1}^{\,n+1}}{2\,\Delta x^{\,2}}\,, \\ j\ =\ 1,\,\ldots,\,N-1\,, \qquad n\ \geq\ 0\,.
\end{multline}
We underline that both schemes \eqref{eq:cnn} and \eqref{eq:ccn} are second order accurate in space \emph{and} in time, \ie the consistency error is $\O(\Delta t^{\,2}\ +\ \Delta x^{\,2})\,$. However, there is a major advantage of the cCN scheme \eqref{eq:ccn} over the classical CN scheme \eqref{eq:cnn} in the fact that cCN is linear with respect to quantities evaluated at the upcoming time layer $t\ =\ t^{\,n+1}$ provided that $k(u)$ is an affine function of $u\,$. This fact can be used to simplify the resolution procedure without destroying the accuracy of the CN scheme. Otherwise, for more general diffusion coefficients $k(u)$ the success of operation depends on the easiness to solve nonlinear equation \eqref{eq:ccn}. It goes without saying that information propagates instantaneously in both CN and cCN schemes.


\subsection{Information propagation speed}

Let us discuss now an important issue of the information propagation speed in the discretized heat equation \eqref{eq:heat1d}. As the initial condition we take the following grid function:
\begin{equation*}
  u_{\,j}^{\,0}\ =\ \begin{dcases}
    \ 1\,,& \qquad j\ =\ 0\,, \\
    \ 0\,,& \qquad j\ \neq\ 0\,,
  \end{dcases}
\end{equation*}
which corresponds to the discrete \textsc{Dirac} function. In all fully implicit schemes (such as \eqref{eq:imp} and \eqref{eq:cn}) the grid function $\bigl\{u_{\,j}^{\,1}\bigr\}_{\,j\,=\,0}^{N}$ will generally have non-zero values in all nodes (modulo perhaps homogeneous boundary conditions). Thus, we can conclude that information has spreaded instantaneously. On the other hand, as it is illustrated in Figure~\ref{fig:grid} with grey and red circles, in explicit discretizations the information propagates one cell to the left and one cell to the right in one time step. Thus, its speed $c_s$ can be estimated as
\begin{equation*}
  c_s\ =\ \frac{\Delta x}{\Delta t}\ \underbrace{\geq}_{\mathrm{CFL}}\ \frac{2\,\nu}{\Delta x}\,.
\end{equation*}
Thus, the value of $c_s$ is finite. Of course, in the limit $\Delta x\ \to\ 0$ we recover the infinite propagation speed, but let us not forget that computations are always run for a \emph{finite} value of $\Delta x\,$.

We arrived to an interesting conclusion. Even if the continuous heat equation \eqref{eq:heat} possesses an unphysical property, it can be corrected if we use a judicious (in this case \emph{explicit}) discretization. This is the main reason why we privilege explicit schemes in time. However, these schemes are subject to severe stability restrictions. The rest of the manuscript is devoted to the question how to overcome the stability limit?

There is another computational advantage of explicit schemes over the implicit ones. Namely, explicit methods can be easily parallelized and they allow to achieve almost perfect scaling on HPC systems \cite{Chetverushkin2012}. Indeed, the computational domain can be split into sub-domains, each sub-domain being handled by a separate processor. Since the stencil is local, only direct neighbours are involved in individual computations. The communication among various processes is almost minimal since only boundary nodes have to be shared. This is another good reason to privilege explicit schemes over the implicit ones.


\section{Improved explicit schemes}
\label{sec:tricks}

Below we present several alternative methods which were specifically designed to overcome the stability limitation of the standard explicit scheme \eqref{eq:exp}.


\subsection{Dufort--Frankel method}
\label{sec:dufort}

\begin{figure}
  \centering
  \includegraphics[width=0.59\textwidth]{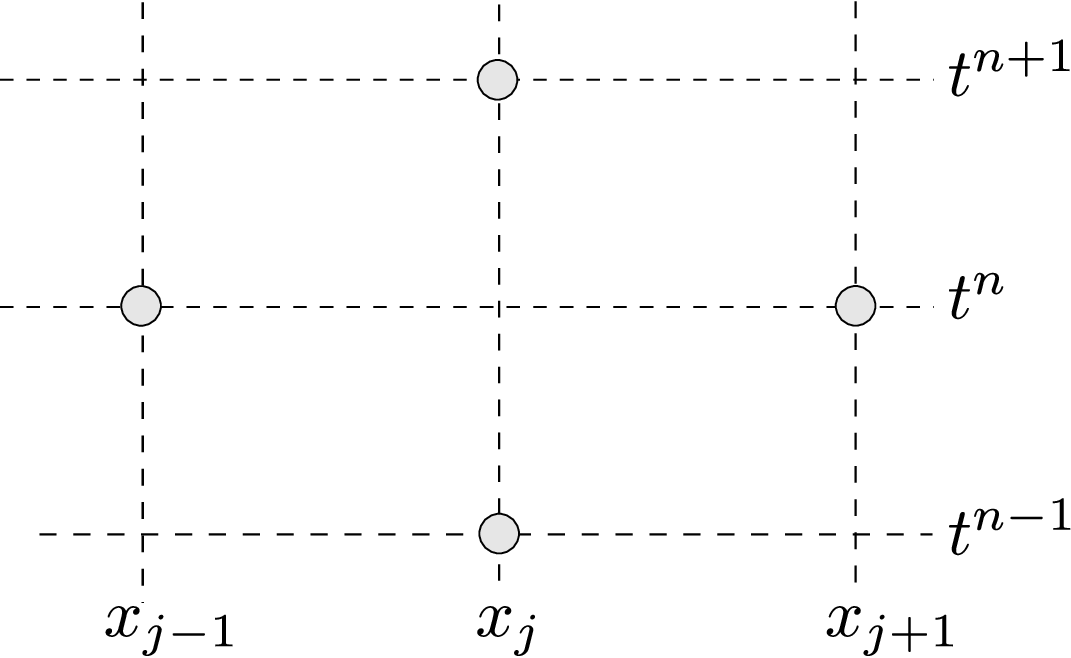}
  \caption{\small\em Stencil of the \textsc{Dufort}--\textsc{Frankel} \eqref{eq:dufort} finite difference scheme.}
  \label{fig:dufort}
\end{figure}

Let us take the unconditionally unstable leap-frog scheme \eqref{eq:leap} and slightly modify it to obtain the so-called \textsc{Dufort}--\textsc{Frankel} method:
\begin{equation}\label{eq:dufort}
  \frac{u_{\,j}^{\,n+1}\ -\ u_{\,j}^{\,n-1}}{2\,\Delta t}\ =\ \nu\;\frac{u_{\,j-1}^{\,n}\ -\ \bigl(u_{\,j}^{\,n-1}\ +\ u_{\,j}^{\,n+1}\bigr)\ +\ u_{\,j+1}^{\,n}}{\Delta x^{\,2}}\,, \qquad j\ =\ 1,\,\ldots,\,N-1\,, \qquad n\ \geq\ 0\,,
\end{equation}
where we made a replacement
\begin{equation*}
  2\,u_{\,j}^{\,n}\ \hookleftarrow\ u_{\,j}^{\,n-1}\ +\ u_{\,j}^{\,n+1}\,.
\end{equation*}
The scheme \eqref{eq:dufort} has the stencil depicted in Figure~\ref{fig:dufort}. At the first glance the scheme \eqref{eq:dufort} looks like an implicit scheme, however, it is not truly the case. Equation \eqref{eq:dufort} can be easily solved for $u_{\,j}^{\,n+1}$ to give the following discrete dynamical system:
\begin{equation*}
  u_{\,j}^{\,n+1}\ =\ \frac{1\ -\ \lambda}{1\ +\ \lambda}\;u_{\,j}^{\,n-1}\ +\ \frac{\lambda}{1\ +\ \lambda}\;\bigl(u_{\,j+1}^{\,n}\ +\ u_{\,j-1}^{\,n}\bigr)\,,
\end{equation*}
where
\begin{equation*}
  \lambda\ \eqdef\ 2\,\nu\;\frac{\Delta t}{\Delta x^{\,2}}\,.
\end{equation*}
The standard \textsc{von Neumann} stability analysis shows that the \textsc{Dufort}--\textsc{Frankel} scheme is \emph{unconditionally stable}.

The consistency error analysis of the scheme \eqref{eq:dufort} shows the following interesting result:
\begin{multline*}
  \L_{\,j}^{\,n}\ =\ \underbrace{\nu\;\frac{\Delta t^{\,2}}{\Delta x^{\,2}}}_{\equiv\ \tau}\;u_{\,t\,t}\ +\ \underbrace{u_{\,t}\ -\ \nu\,u_{\,x\,x}}_{\eqref{eq:heat1d}}\ +\\
  \sixth\,\Delta t^{\,2}\,u_{\,t\,t\,t}\ -\ \twelwth\,\nu\,\Delta x^{\,2}\,u_{\,x\,x\,x\,x}\ -\ \twelwth\,\nu\,\Delta t^{\,2}\,\Delta x\,u_{\,x\,x\,x\,t\,t}\ +\ \O\Bigl(\frac{\Delta t^{\,4}}{\Delta x^{\,2}}\Bigr)\,,
\end{multline*}
where
\begin{equation*}
  \L_{\,j}^{\,n}\ \eqdef\ \frac{u_{\,j}^{\,n+1}\ -\ u_{\,j}^{\,n-1}}{2\,\Delta t}\ -\ \nu\;\frac{u_{\,j-1}^{\,n}\ -\ \bigl(u_{\,j}^{\,n-1}\ +\ u_{\,j}^{\,n+1}\bigr)\ +\ u_{\,j+1}^{\,n}}{\Delta x^{\,2}}\,.
\end{equation*}
So, from the asymptotic expansion for $\L_{\,j}^{\,n}$ we obtain that the \textsc{Dufort}--\textsc{Frankel} scheme is second order accurate in time \emph{and}
\begin{itemize}
  \item First order accurate in space if $\Delta t\ \propto\ \Delta x^{\,3/2}$
  \item Second order accurate in space if $\Delta t\ \propto\ \Delta x^{\,2}$
\end{itemize}
However, the \textsc{Dufort}--\textsc{Frankel} scheme is unconditionally consistent with the so-called \emph{hyperbolic heat conduction equation}:
\begin{equation*}
  \tau\,u_{\,t\,t}\ +\ u_{\,t}\ -\ \nu\,u_{\,x\,x}\ =\ 0\,.
\end{equation*}
We shall return to this equation below. At this stage we only mention that information propagates with the finite speed in hyperbolic models.

\subsection{Saulyev method}
\label{sec:saulyev}

In this Section we describe a not so widely known method proposed by \textsc{Saulyev} in \cite{Saulyev1960} to integrate parabolic equations. For simplicity we focus on the $1-$dimensional heat equation \eqref{eq:heat1d}. The first idea of this method consists in rewriting the discrete second (spatial) derivative as
\begin{equation*}
  u_{\,x\,x}\bigr\vert_{x\,=\,x_j}\ \approx\ \frac{u_{\,j+1}\ -\ 2\,u_{\,j}\ +\ u_{\,j-1}}{\Delta x^{\,2}}\ \equiv\ \frac{\dfrac{u_{\,j+1}\ -\ u_{\,j}}{\Delta x}\ -\ \dfrac{u_{\,j}\ -\ u_{\,j-1}}{\Delta x}}{\Delta x}\,.
\end{equation*}
In the finite difference formula above we do not specify intentionally the time layer number. The next trick consists in writing the following asymmetric finite difference approximation:
\begin{equation*}
  \frac{u_{\,j}^{\,n+1}\ -\ u_{\,j}^{\,n}}{\Delta t}\ =\ \nu\;\frac{\dfrac{u_{\,j+1}^{\,n}\ -\ u_{\,j}^{\,n}}{\Delta x}\ -\ \dfrac{u_{\,j}^{\,n+1}\ -\ u_{\,j-1}^{\,n+1}}{\Delta x}}{\Delta x}\,,
\end{equation*}
or after simplifications we obtain
\begin{equation}\label{eq:s1}
  \frac{u_{\,j}^{\,n+1}\ -\ u_{\,j}^{\,n}}{\Delta t}\ =\ \nu\;\frac{u_{\,j+1}^{\,n}\ -\ \bigl(u_{\,j}^{\,n}\ +\ u_{j}^{\,n+1}\bigr)\ +\ u_{\,j-1}^{\,n+1}}{\Delta x^{\,2}}\,.
\end{equation}
The last difference relation is slightly different from the \textsc{Dufort}--\textsc{Frankel} method presented above. Moreover, the relation written above is not consistent with the original equation \eqref{eq:heat1d}. That is why we consider the next time layer $t\ =\ t^{\,n+2}$ and we apply symmetrically the same tricks, \ie
\begin{equation*}
  \frac{u_{\,j}^{\,n+2}\ -\ u_{\,j}^{\,n+1}}{\Delta t}\ =\ \nu\;\frac{\dfrac{u_{\,j+1}^{\,n+2}\ -\ u_{\,j}^{\,n+2}}{\Delta x}\ -\ \dfrac{u_{\,j}^{\,n+1}\ -\ u_{\,j-1}^{\,n+1}}{\Delta x}}{\Delta x}\,,
\end{equation*}
and after simplifications we have
\begin{equation}\label{eq:s2}
  \frac{u_{\,j}^{\,n+2}\ -\ u_{\,j}^{\,n+1}}{\Delta t}\ =\ \nu\;\frac{u_{\,j+1}^{\,n+2}\ -\ \bigl(u_{\,j}^{\,n+2}\ +\ u_{j}^{\,n+1}\bigr)\ +\ u_{\,j-1}^{\,n+1}}{\Delta x^{\,2}}\,.
\end{equation}
Without any surprise the relation \eqref{eq:s2} does not approximate equation \eqref{eq:heat1d} either. However, both relations \eqref{eq:s1} and \eqref{eq:s2} constitute the so-called \textsc{Saulyev} scheme and are called the first and second stages of \textsc{Saulyev} method correspondingly.

In order to perform the approximation error analysis we take the sum of \eqref{eq:s1}, \eqref{eq:s2} and we divide it by two:
\begin{equation*}
  \L_{\,j}^{\,n}\ =\ \frac{u_{\,j}^{\,n+2}\ -\ u_{\,j}^{\,n}}{2\,\Delta t}\ =\ \nu\;\frac{u_{\,j}^{\,n}\ -\ u_{\,j+1}^{\,n}\ -\ 2\,u_{\,j-1}^{\,n+1}\ +\ 2\,u_{\,j}^{\,n+1}\ +\ u_{\,j}^{\,n+2}\ -\ u_{\,j+1}^{\,n+2}}{2\,\Delta x^{\,2}}\,.
\end{equation*}
After applying local \textsc{Taylor} expansions we obtain
\begin{multline*}
  \L_{\,j}^{\,n}\ =\ u_{\,t}\ -\ \nu\,u_{\,x\,x}\ -\ \twelwth\,\nu\,\Delta x^{\,2}\,u_{\,x\,x\,x\,x}\ +\\ 
  \bigl[\,\twothird\,u_{\,t\,t\,t}\ -\ \threefourth\,\nu\,u_{\,x\,x\,t\,t}\,\bigr]\,\Delta t^{\,2}\ -\ \half\,\nu\;\frac{\Delta t^{\,2}}{\Delta x}\;u_{\,x\,t\,t}\ +\ \O(\Delta t\,\Delta x^{\,2}\ +\ \Delta t^{\,2}\,\Delta x)\,.
\end{multline*}
From the last asymptotic expansion we arrive at an important result: \textsc{Saulyev} scheme is second order accurate in space if $\ \Delta t\ =\ \O(\Delta x^{\,3/2})\,$. We underline the fact that this condition coming from the accuracy requirements is weaker than the usual CFL restriction \eqref{eq:cfl}. For instance, if the user is ready to sacrifice the spatial accuracy to the first order $\O(\Delta x)\,$, then it is sufficient to take $\ \Delta t\ =\ \O(\Delta x)\,$.

\begin{figure}
  \centering
  \includegraphics[width=0.69\textwidth]{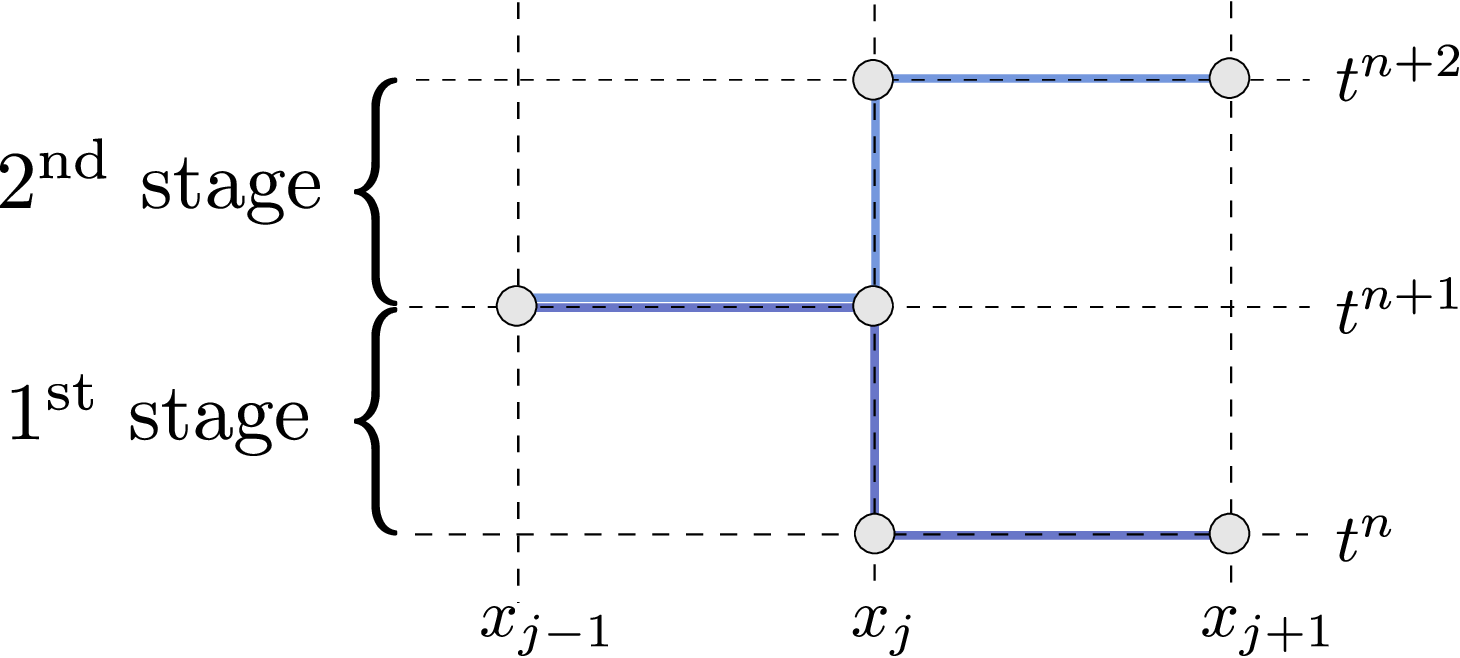}
  \caption{\small\em Stencil of the \textsc{Saulyev} finite difference scheme, which consists of two stages \eqref{eq:s1} and \eqref{eq:s2}.}
  \label{fig:saul}
\end{figure}

Without proof we report that \textsc{Saulyev}'s scheme is unconditionally stable (see \cite{Saulyev1960} for more details). The stencil of \textsc{Saulyev}'s scheme is depicted in Figure~\ref{fig:saul}.

\subsubsection{Resolution procedure}

At the first glance \textsc{Saulyev} scheme appears as an implicit one since each relation \eqref{eq:s1} and \eqref{eq:s2} contains two terms from the following time layer ($t\ =\ t^{\,n+1}$ and $t\ =\ t^{\,n+2}$ correspondingly). However, this scheme can be recast in an almost explicit form using judicious recurrence relations.

Consider the first stage \eqref{eq:s1} of \textsc{Saulyev}'s scheme. Similarly to Section~\ref{sec:dufort} we introduce for simplicity the parameter
\begin{equation*}
  \lambda\ \eqdef\ \nu\;\frac{\Delta t}{\Delta x^{\,2}}\,.
\end{equation*}
From difference relation \eqref{eq:s1} we find
\begin{equation}\label{eq:rec1}
  u_{\,j}^{\,n+1}\ =\ \frac{1\ -\ \lambda}{1\ +\ \lambda}\;u_{\,j}^{\,n}\ +\ \frac{\lambda}{1\ +\ \lambda}\;u_{\,j+1}^{\,n}\ +\ \frac{\lambda}{1\ +\ \lambda}\;u_{\,j-1}^{\,n+1}\,.
\end{equation}
The first stage of \textsc{Saulyev}'s scheme is computed in rightwards direction (increasing index $j\nearrow$). From the left boundary condition we find first the value
\begin{equation*}
  u_{\,0}^{\,n+1}\ =\ \psi_{\,\mathrm{l}}\,(t^{\,n+1})\,.
\end{equation*}
This allows us to compute $u_{\,1}^{\,n+1}\,$, $u_{\,2}^{\,n+1}\,$, \etc thanks to the recurrence relation \eqref{eq:rec1}. At the final step, the value $u_{\,N}^{\,n+1}\,$ is computed directly from the right boundary condition:
\begin{equation*}
  u_{\,N}^{\,n+1}\ =\ \psi_{\,\mathrm{r}}\,(t^{\,n+1})\,.
\end{equation*}
This completes the first stage of computations.

\begin{remark}
For some types of boundary conditions (\eg \textsc{Robin}-type) \textsc{Saulyev}'s scheme might require solution of a small dimensional (typically $2\,\times\,2$) linear system of algebraic equations to initiate the recurrence \eqref{eq:rec1}.
\end{remark}

Let us make explicit now the second stage \eqref{eq:s2} of \textsc{Saulyev}'s method. For this purpose we solve relation \eqref{eq:s2} with respect to $u_{\,j}^{\,n+2}\,$:
\begin{equation}\label{eq:rec2}
  u_{\,j}^{\,n+2}\ =\ \frac{1\ -\ \lambda}{1\ +\ \lambda}\;u_{\,j}^{\,n+1}\ +\ \frac{\lambda}{1\ +\ \lambda}\;u_{\,j-1}^{\,n+1}\ +\ \frac{\lambda}{1\ +\ \lambda}\;u_{\,j+1}^{\,n+2}\,.
\end{equation}
Now it is getting clear that during the second stage of \textsc{Saulyev}'s scheme we proceed in the leftwards direction (decreasing $j\searrow$). From the right boundary condition we find first
\begin{equation*}
  u_{\,N}^{\,n+2}\ =\ \psi_{\,\mathrm{r}}\,(t^{\,n+2})\,.
\end{equation*}
It allows us to compute $u_{\,N-1}^{\,n+2}\,$, $u_{\,N-2}^{\,n+2}\,$, \etc using the recurrence relation \eqref{eq:rec2}. At the final step, the value $u_{\,0}^{\,n+2}\,$ is computed directly from the left boundary condition:
\begin{equation*}
  u_{\,0}^{\,n+2}\ =\ \psi_{\,\mathrm{l}}\,(t^{\,n+2})\,.
\end{equation*}
As a result we obtain a fully explicit resolution scheme without stability related limitations. We notice however that \textsc{Saulyev}'s scheme provides consistent results only every second time step or, in other words, after the successive completion of both stages \eqref{eq:rec1} and \eqref{eq:rec2}. The intermediate result is not consistent with the equation \eqref{eq:heat1d}.


\subsection{Hyperbolization method}
\label{sec:hyper}

We saw above that the \textsc{Dufort}--\textsc{Frankel} scheme is a hidden way to add a small amount of `hyperbolicity' into the model \eqref{eq:heat1d}. In this Section we shall invert the order of operations: first, we perturb the equation \eqref{eq:heat1d} in an ad-hoc way and only after we discretize it with a suitable method.

Consider the $1-$dimensional heat equation \eqref{eq:heat1d} that we are going to perturb by adding a small term containing the second derivative in time:
\begin{equation}\label{eq:hyp}
  \tau\,u_{\,t\,t}\ +\ u_{\,t}\ -\ \nu\,u_{\,x\,x}\ =\ 0\,.
\end{equation}
This is the \emph{hyperbolic heat equation} already familiar to us since it appeared in the consistency analysis of the \textsc{Dufort}--\textsc{Frankel} scheme. Here we perform a singular perturbation by assuming that
\begin{equation*}
  \norm{\tau\,u_{\,t\,t}}\ \ll\ \norm{u_{\,t}}\,.
\end{equation*}
The last condition physically means that the new term has only limited influence on the solution of equation \eqref{eq:hyp}. Here $\tau$ is a small ad-hoc parameter whose value is in general related to physical and discretization parameters $\tau\ =\ \tau\,(\nu,\,\Delta x,\,\Delta t)\,$.

\begin{remark}
One can notice that equation \eqref{eq:hyp} is second order in time, thus, it requires two initial conditions to obtain a well-posed initial value problem. However, the parabolic equation \eqref{eq:heat1d} is only first order in time and it only requires the knowledge of the initial temperature field distribution. When we solve the hyperbolic equation \eqref{eq:hyp}, the missing initial condition is simply chosen to be
\begin{equation*}
  u_{\,t}\,\bigr\vert_{t\,=\,0}\ =\ 0\,.
\end{equation*}
\end{remark}


\subsubsection{Dispersion relation analysis}

The classical dispersion relation analysis looks at plane wave solutions:
\begin{equation}\label{eq:ans}
  u(x,\,t)\ =\ u_{\,0}\,\ue^{\,\ui\,(\kappa\,x\ -\ \omega\,t)}\,.
\end{equation}
By substituting this solution ansatz into equation \eqref{eq:heat1d} we obtain the following relation between wave frequency $\omega$ and wavenumber $k\,$:
\begin{equation}\label{eq:disp}
  \omega(\kappa)\ =\ -\,\ui\,\nu\,\kappa^{\,2}\,.
\end{equation}
The last relation is called the \emph{dispersion relation} even if the heat equation \eqref{eq:heat1d} is not dispersive but dissipative. The real part of $\omega$ contains information about wave propagation properties (dispersive if $\ \frac{\Re\omega(\kappa)}{\kappa}\ \neq\ \const$ and non-dispersive otherwise) while the imaginary part describes how different modes $\kappa$ dissipate (if $\ \Im\omega(\kappa)\ <\ 0$) or grow (if $\ \Im\omega(\kappa)\ >\ 0\,$). The dispersion relation \eqref{eq:disp} gives the damping rate of different modes.

The same plane wave ansatz \eqref{eq:ans} can be substituted into the hyperbolic heat equation \eqref{eq:hyp} as well to give the following \emph{implicit} relation for the wave frequency $\omega\,$:
\begin{equation*}
  -\,\tau\,\omega^{\,2}\ -\ \ui\,\omega\ +\ \nu\,\kappa^{\,2}\ =\ 0\,.
\end{equation*}
By solving this quadratic equation with complex coefficients for $\omega\,$, we obtain two branches:
\begin{equation*}
  \omega_{\,\pm}\,(\kappa)\ =\ \frac{-\,\ui\ \pm\ \sqrt{4\,\nu\,\kappa^{\,2}\,\tau\ -\ 1}}{2\,\tau}\,.
\end{equation*}
This dispersion relation will be analyzed asymptotically with $\tau\ \ll\ 1$ being the small parameter. The branch $\omega_{\,-}\,(\kappa)$ is not of much interest to us since it is constantly damped, \ie
\begin{equation*}
  \omega_{\,-}\,(\kappa)\ =\ -\;\frac{\ui}{\tau}\ +\ \O(1)\,.
\end{equation*}
It is much more instructive to look at the positive branch $\omega_{\,+}\,(\kappa)\,$:
\begin{equation*}
  \omega_{\,+}\,(\kappa)\ =\ -\,\ui\,\nu\,\kappa^{\,2}\,\bigl[\,1\ +\ \nu\,\kappa^{\,2}\,\tau\ +\ 2\,\nu^{\,2}\,\kappa^{\,4}\,\tau^{\,2}\ +\ \O(\tau^{\,3})\,\bigr]\,.
\end{equation*}
The last asymptotic expansion shows that for small values of parameter $\tau$ we obtain a valid asymptotic approximation of the dispersion relation \eqref{eq:disp} for the heat equation \eqref{eq:heat1d}.


\subsubsection{Discretization}

Equation \eqref{eq:hyp} will be discretized on the same stencil as the leap-frog scheme \eqref{eq:leap} (see Figure~\ref{fig:leap}):
\begin{multline}\label{eq:hyps}
  \L_{\,j}^{\,n}\ \eqdef\ \tau\;\frac{u_{\,j}^{\,n+1}\ -\ 2\,u_{\,j}^{\,n}\ +\ u_{\,j}^{\,n-1}}{\Delta t^{\,2}}\ +\ \frac{u_{\,j}^{\,n+1}\ -\ u_{\,j}^{\,n-1}}{2\,\Delta t}\ -\ \nu\;\frac{u_{\,j+1}^{\,n}\ -\ 2\,u_{\,j}^{\,n}\ +\ u_{\,j-1}^{\,n}}{\Delta x^{\,2}}\ =\ 0\,,\\ \qquad j\ =\ 1,\,\ldots,\,N-1\,, \qquad n\ \geq\ 0\,,
\end{multline}
The last scheme is consistent with hyperbolic heat equation \eqref{eq:hyp} to the second order in space and in time $\O(\Delta t^{\,2}\ +\ \Delta x^{\,2})\,$. Indeed, using the standard \textsc{Taylor} expansions we obtain
\begin{multline*}
  \L_{\,j}^{\,n}\ =\ \tau\,u_{\,t\,t}\ +\ u_{\,t}\ -\ \nu\,u_{\,x\,x}\\
  -\ \frac{\nu}{12}\;\Delta x^{\,2}\,u_{\,x\,x\,x\,x}\ +\ \Delta t^{\,2}\,\Bigl[\,\sixth\,u_{\,t\,t\,t}\ +\ \twelwth\,\tau\,u_{\,t\,t\,t\,t}\,\Bigr]\ +\ \O(\Delta t^{\,4}\ +\ \Delta x^{\,4})\,.
\end{multline*}

The stability of the scheme \eqref{eq:hyps} was studied in \cite{Chetverushkin2012} and the following stability condition was obtained:
\begin{equation*}
  \frac{\Delta t}{\Delta x}\ \leq\ \sqrt{\frac{\tau}{\nu}}\,.
\end{equation*}
By taking, for example, $\tau\ =\ \nu\,\Delta x$ we obtain the following stability condition
\begin{equation*}
  \Delta t\ \leq\ \Delta x^{\,\frac{3}{2}}\,,
\end{equation*}
which is still weaker than the standard parabolic condition \eqref{eq:cfl}. However, it was reported in \cite{Chetverushkin2010, Chetverushkin2015} that stable computations (even in 3D) can be performed even with $\Delta t\ =\ \O(\Delta x)\,$. The authors explain informally this experimental observation by the fact that usual stability conditions are too `pessimistic'.

\begin{remark}
The ad-hoc parameter $\tau$ can be chosen in other ways as well. One popular choice consists in taking
\begin{equation*}
  \tau\ =\ \frac{\Delta x}{c_s}\,,
\end{equation*}
where $c_s$ is the real physical information speed.
\end{remark}


\subsubsection{Error estimate}

It is legitimate to ask the question how far are solutions $u_{\,h}\,(x,\,t)$ to the hyperbolic equation \eqref{eq:hyp} from the solutions $u_{\,p}\,(x,\,t)$ of the parabolic heat equation \eqref{eq:heat} (for the same initial condition). This question for the initial value problem was studied in \cite{Myshetskaya2015} and we shall provide here only the obtained error estimate. Let us introduce the difference between two solutions:
\begin{equation*}
  \delta u\,(x,\,t)\ \eqdef\ u_{\,h}\,(x,\,t)\ -\ u_{\,p}\,(x,\,t)\,.
\end{equation*}
Then, the following estimate holds
\begin{equation*}
  \abs{\delta u\,(x,\,t)}\ \leq\ \tau\,\M\,\Bigl(1\ +\ \frac{2}{\sqrt{\pi}}\Bigr)\,\biggl(8\,\sqrt{2}\,\tau\ +\ \frac{\sqrt[4]{2\,\pi^2}}{2}\;T\biggr)\,,
\end{equation*}
where $T\ >\ 0$ is the time horizon and
\begin{equation*}
  \M\ \eqdef\ \sup_{\Omega_{\,\xi,\,\zeta}}\;\Bigl\vert\,\pd{^{\,2}u_{\,p}}{\,t^{\,2}}(\xi,\,\zeta)\,\Bigr\vert\,,
\end{equation*}
and the domain $\Omega_{\,\xi,\,\zeta}$ is defined as
\begin{equation*}
  \Omega_{\,\xi,\,\zeta}\ \eqdef\ \Bigl\{\,(\xi,\,\zeta)\ :\ 0\ \leq\ \zeta\ \leq\ t\,, \quad
  x\ -\ \frac{t\ -\ \zeta}{\sqrt{\tau}}\ \leq\ \xi\ \leq\ x\ +\ \frac{t\ -\ \zeta}{\sqrt{\tau}}\,\Bigr\}\,.
\end{equation*}


\section{Discussion}
\label{sec:disc}

We saw above that only explicit schemes allow to have the finite speed of information propagation in the discretized version of the heat equation \eqref{eq:heat}. The ease of parallelization along with excellent scaling properties of the codes obtained with explicit schemes constitute another important advantage to privilege explicit schemes over implicit ones \cite{Chetverushkin2012}. However, explicit schemes for parabolic equations suffer from very stringent CFL-type conditions $\Delta t\ =\ \O(\Delta x^{\,2})$ on the time step. That is why it is almost impossible to perform long time simulations (required in \eg building physics applications) using fully explicit schemes such as \eqref{eq:exp}. In order to keep explicit discretizations and overcome stringent CFL-type restrictions, a certain number of \emph{hybrid} schemes were described. These hybrid schemes are based on different ideas. Some schemes rely on the information about the numerical solution on following time layers while keeping the overall scheme explicit using various tricks. Some hybrid schemes (\eg \textsc{Dufort}--\textsc{Frankel} method described in Section~\ref{sec:dufort} and \textsc{Saulyev}'s scheme from Section~\ref{sec:saulyev}) can be even unconditionally stable.

Using the local error analysis we noticed that CFL-improved schemes gain the stability by introducing some weak hyperbolicity into the model \eqref{eq:heat1d}. In particular, it is the case of the \textsc{Dufort}--\textsc{Frankel} scheme, while \textsc{Saulyev} method seems to be rather dispersive. This observation suggests that a new method can be derived by introducing this hyperbolicity in a controlled manner and to discretize the perturbed equation later. The method of hyperbolization deforms the equation operator to achieve desired properties\footnote{By desired properties we mean the finite speed of information propagation along with less stringent CFL-type restrictions in explicit finite difference discretizations.} of the numerical solution. However, in all cases the gain in stability results in some loss of accuracy in representing the original continuous operator \eqref{eq:heat}. The trade-off between the stability and accuracy has to be made by the end user. The second order accuracy is equivalent to the classical parabolic CFL-type condition. However, the user can choose to degrade intentionally the accuracy to relax the stability restriction up to hyperbolic-type conditions $\Delta t\ =\ \O(\Delta x)$ (originally found in \cite{Courant1928}) and even beyond. It is not difficult to see that hybrid schemes presented in this study can be easily generalized to nonlinear cases with source terms in one and more spatial dimensions.

The main goal of this manuscript was to communicate and attract community's attention to these improved discretizations, which can be used in modern building physics simulations where typical time scales are measured rather in months or even years. Our preference goes perhaps to the method of hyperbolization since it has been successfully applied (and validated) even to compressible \textsc{Navier}--\textsc{Stokes} and MHD equations \cite{Chetverushkin2015}. Another advantage of hyperbolization technique is that it can be mathematically derived for gas dynamics from the kinetic theory of \textsc{Boltzmann}.

While preparing this text, the Author discovered also another strategy to relax the the time step limits. This idea can be traced back to the scientific school of R.~\textsc{Temam} who proposed to separate the scales for the purposes of numerical simulations. Then, high frequencies, which impose severe stability limits, are treated separately. For a modern introduction to these approaches we refer to \cite{Brachet2016}. There exist also explicit time integration methods\textsc{The so-called \textsc{Runge}--\textsc{Kutta}--\textsc{Tchebyshev} methods.} especially designed to have extended stability limits. This research direction was pioneered also by V.~K.~\textsc{Saulyev} in 1960 \cite{Saulyev1960}.


\subsection*{Acknowledgments}
\addcontentsline{toc}{subsection}{Acknowledgments}

The author would like to thank Professor B.N.~\textsc{Chetverushkin} whose talks available on-line greatly inspired the preparation of this manuscript. I would like also to thank my collaborators J.~\textsc{Berger}, S.~\textsc{Gasparin} and N.~\textsc{Mendes} for bringing my attention to heat and mass conduction problems.
\bigskip


\addcontentsline{toc}{section}{References}
\bibliographystyle{abbrv}
\bibliography{biblio}

\end{document}